\documentclass[11pt,oneside]{article}

\usepackage[T1]{fontenc}
\usepackage[margin=1in,includefoot]{geometry}
\usepackage{mathpazo}
\usepackage{xcolor}

\newcommand{\red}[1]{\textcolor{red}{#1}}

\usepackage[hidelinks]{hyperref}

\usepackage{caption}
\usepackage{graphicx}

\let\oldincludegraphics\includegraphics%
\renewcommand{\includegraphics}[2][]{\IfFileExists{#2}{\oldincludegraphics[#1]{#2}}{\red{[FILE NOT FOUND]}}}

\usepackage{amsmath}
\usepackage{amssymb}
\usepackage{amsthm}

\DeclareMathOperator{\E}{E}

\DeclareMathOperator{\Cov}{Cov}

\DeclareMathOperator{\Var}{Var}

\newcommand{\abs}[1]{\lvert#1\rvert}
\newcommand{\Dcal}{\mathcal{D}}
\newcommand{\der}{\mathrm{d}}

\newcommand{\Ncal}{\mathcal{N}}
\newcommand{\R}{\mathbb{R}}
\newcommand{\seps}{\sigma_\epsilon}
\newcommand{\sm}{\sigma}
\newcommand{\sz}{\sigma_0}
\newcommand{\veps}{\seps^2}
\newcommand{\vm}{\sm^2}
\newcommand{\vz}{\sz^2}
\renewcommand{\epsilon}{\varepsilon}

\let\oldleft\left
\let\oldright\right
\renewcommand{\left}{\mathopen{}\mathclose\bgroup\oldleft}
\renewcommand{\right}{\aftergroup\egroup\oldright}

\newtheorem{corollary}{Corollary}
\newtheorem{lemma}{Lemma}

\newtheorem{theorem}{Theorem}

\newif\ifbodyproofs
\bodyproofstrue

\usepackage{natbib}
\bodyproofsfalse

\title{Estimating sample paths of Gauss-Markov\par processes from noisy data}
\author{%
Benjamin Davies\thanks{
Department of Economics, Stanford University; bldavies@stanford.edu.
}
}
\date{Draft version: \today}

\begin{document}

\maketitle

\begin{abstract}
    \noindent
    I derive the pointwise conditional means and variances of an arbitrary Gauss-Markov process, given noisy observations of points on a sample path.
    These moments depend on the process's mean and covariance functions, and on the conditional moments of the sampled points.
    I study the Brownian motion and bridge as special cases.
    
    \vskip\baselineskip
\end{abstract}

\section{Introduction}

Suppose we observe data~$\Dcal=\{(x_i,y_i)\}_{i=1}^n$ generated by the process
\begin{equation}
    y_i=f(x_i)+\epsilon_i, \label{eq:dgp}
\end{equation}
where~$x_i\ge0$ is non-decreasing in~$i$, where~$f:[0,\infty)\to\R$ is unknown, and where the errors~$\epsilon_i$ are jointly normally distributed (hereafter ``Gaussian'') with~$\E[\epsilon_i\mid x_1,x_2\ldots,x_n]=0$ independently of~$\{f(x)\}_{x\ge0}$.
We use~$\Dcal$ to construct pointwise estimates
\[ \hat{f}_\Dcal(x)\equiv\E[f(x)\mid\Dcal] \]
with mean squared error (MSE)
\begin{align}
    \E\left[\left(f(x)-\hat{f}_\Dcal(x)\right)^2\mid\Dcal\right]
    &= \Var(f(x)\mid\Dcal). \label{eq:mse}
\end{align}
In this note, I derive expressions for~$\E[f(x)\mid\Dcal]$ and~$\Var(f(x)\mid\Dcal)$ when~$\{f(x)\}_{x\ge0}$ is a sample path of a Gauss-Markov process.
Such processes have two defining properties:
\begin{enumerate}

    \item[(G)]
    For every finite subset~$X\subset[0,\infty)$, the vector~$(f(x))_{x\in X}\in\R^{\abs{X}}$ is multivariate Gaussian;

    \item[(M)]
    If~$x<x'<x''$, then~$f(x)$ and~$f(x'')$ are conditionally independent given~$f(x')$.

\end{enumerate}
Property~(G) implies that~$f(x)\mid\Dcal$ is Gaussian, and so its distribution is fully determined by its mean~$\E[f(x)\mid\Dcal]$ and variance~$\Var(f(x)\mid\Dcal)$.
Theorem~\ref{thm:moments-arbitrary} expresses these moments in terms of the means and (co)variances of the~$f(x_i)\mid\Dcal$.
This allows me to construct the estimate~$\hat{f}_\Dcal(x)$ and its MSE at all points~$x\ge0$.
This estimate optimally extrapolates from, or interpolates between, the observations~$(x_i,y_i)$ in~$\Dcal$.%
\footnote{
The estimate is ``optimal'' in that it minimizes the MSE~\eqref{eq:mse} for all~$x\ge0$.
}

I let these observations be noisy, with Gaussian errors~$\epsilon_i=y_i-f(x_i)$.
This allows me to extend analyses that assume observations have no noise \cite[e.g.,][]{Bardhi-2024-ECTA,Callander-2011-AER,Carnehl-Schneider-2023-}.%
\footnote{
\citet[][p.\! 16]{Rasmussen-Williams-2006-} derive expressions for~$\E[f(x)\mid\Dcal]$ and~$\Var(f(x)\mid\Dcal)$ when~$\{f(x)\}_{x\ge0}$ follows a Gaussian (but not necessarily Markov) process and the errors~$\epsilon_i$ are iid.
I impose the Markov property~(M) to obtain (relatively) closed-form expressions for~$\E[f(x)\mid\Dcal]$ and~$\Var(f(x)\mid\Dcal)$.
I also allow for arbitrary (co)variances in the~$\epsilon_i$.
}
I also let the sampled points~$x_i$ be less or greater than the target point~$x$.
This contrasts with \cite{Davies-2024-}, who studies sequential learning from noisy observations of a sample path.
Such learning always involves extrapolation, whereas I allow for interpolation.


\section{Preliminaries}
\label{sec:prelims}

The data~$\Dcal=\{(x_i,y_i)\}_{i=1}^n$ contain noisy observations~$y_i=f(x_i)+\epsilon_i$ of the values~$f(x_i)$.
These observations equal the sum of two Gaussian random variables and so are Gaussian too.
Moreover, by property~(G), the vector~$(f(x_1),f(x_2),\ldots,f(x_n),f(x))$ is multivariate Gaussian for all~$x\ge0$.
It follows that
\[ (y_1,y_2,\ldots,y_n,f(x))=(f(x_1),f(x_2),\ldots,f(x_n),f(x))+(\epsilon_1,\epsilon_2,\ldots,\epsilon_n,f(x)) \]
is also multivariate Gaussian.
Consequently, we can construct the conditional distribution of~$f(x)$ given~$y_1,y_2,\ldots,y_n$ using a well-known result about multivariate Gaussian variables:
\begin{lemma}
    \label{lem:conditional-many}
    Let~$n_1\ge1$ and~$n_2\ge1$ be integers, and let~$z\sim\Ncal(\mu,\Sigma)$ be multivariate Gaussian with mean~$\mu\in\R^{n_1+n_2}$ and variance~$\Sigma\in\R^{(n_1+n_2)\times(n_1+n_2)}$.
    Partition~$z=(z_1,z_2)$ into vectors~$z_1\in\R^{n_1}$ and~$z_2\in\R^{n_2}$, and let~$\mu=(\mu_1,\mu_2)$ and
    \[ \Sigma=\begin{pmatrix}
        \Sigma_{11} & \Sigma_{12} \\
        \Sigma_{21} & \Sigma_{22}
    \end{pmatrix} \]
    be the corresponding partitions of~$\mu$ and~$\Sigma$.
    If~$\Sigma_{22}$ is invertible, then
    \begin{equation}
        z_1\mid z_2\sim\Ncal\left(\mu_1+\Sigma_{12}\Sigma_{22}^{-1}(z_2-\mu_2),\ \Sigma_{11}-\Sigma_{12}\Sigma_{22}^{-1}\Sigma_{21}\right). \label{eq:conditional-many}
    \end{equation}
\end{lemma}

See \citet[p.\! 87]{Bishop-2006-} or \citet[p.\! 55]{DeGroot-2004-} for proofs of this lemma, and Appendix~\ref{sec:proofs} for proofs of my other results.

Substituting~$z_1=f(x)$ and~$z_2=(y_1,y_2,\ldots,y_n)$ into~\eqref{eq:conditional-many} provides expressions for the moments of~$f(x)\mid\Dcal=z_1\mid z_2$.
I refine these expressions by imposing properties~(G) and~(M).

Property~(G) comes from~$\{f(x)\}_{x\ge0}$ being a sample path of a Gaussian process.%
\footnote{
Gaussian processes are stochastic processes satisfying by property~(G).
For more information on these processes and their applications, see Section~6.4 of \citet{Bishop-2006-} or Chapter~2 of \cite{Rasmussen-Williams-2006-}.
}
This process can be characterized by
\begin{enumerate}

    \item
    A mean function~$m:[0,\infty)\to\R$ with~$m(x)\equiv\E[f(x)]$ for all~$x\ge0$, and

    \item
    A covariance function~$C:[0,\infty)^2\to\R$ with~$C(x,x')\equiv\Cov(f(x),f(x'))$ for all~$x,x'\ge0$.

\end{enumerate}
For convenience, I define a variance function~$V:[0,\infty)\to\R$ by~$V(x)\equiv C(x,x)$ for all~$x\ge0$.
The values of~$m(x)$, $V(x)$, and~$C(x,x')$ are known for all~$x,x'\ge0$, but the values of~$f(x)$ are not.

Property~(M) comes from~$\{f(x)\}_{x\ge0}$ being a sample path of a Markov process.
It allows me to focus on the conditional distribution of~$f(x)$ given at most two values~$f(x_i)$: those with~$x_i$ closest to~$x$.
Lemma~\ref{lem:conditional-12} characterizes this conditional distribution.

\begin{lemma}
    \label{lem:conditional-12}
    Let~$\{f(x)\}_{x\ge0}$ be a sample path of a Gauss-Markov process, let~$x,x',x''\ge0$ be arbitrary, and define
    \begin{equation}
        w(x,x',x'')\equiv \frac{1}{V(x')V(x'')-C(x',x'')^2}
        \begin{bmatrix} C(x,x')V(x'')-C(x,x'')C(x',x'') \\ C(x,x'')V(x')-C(x,x')C(x',x'') \end{bmatrix}.
    \end{equation}
    Then~$f(x)\mid f(x')$ is Gaussian with mean
    \begin{equation}
        \E[f(x)\mid f(x')]=m(x)+\frac{C(x,x')}{V(x')}(f(x')-m(x')) \label{eq:conditional-mean-1}
    \end{equation}
    and variance
    \begin{equation}
        \Var(f(x)\mid f(x'))=V(x)-\frac{C(x,x')^2}{V(x')}, \label{eq:conditional-variance-1}
    \end{equation}
    and~$f(x)\mid f(x'),f(x'')$ is Gaussian with mean
    \[ \E[f(x)\mid f(x'),f(x')]=m(x)+w(x,x',x'')^T\begin{bmatrix} f(x')-m(x') \\ f(x'')-m(x'') \end{bmatrix} \]
    and variance
    \[ \Var(f(x)\mid f(x'),f(x''))=V(x)-w(x,x',x'')^T\begin{bmatrix} C(x,x') \\ C(x,x'') \end{bmatrix}. \]
\end{lemma}
\ifbodyproofs\begin{proof}[Proof of Lemma~\ref{lem:conditional-12}]
    Now~$f(x)$, $f(x')$, and~$f(x'')$ are jointly Gaussian by property~(G).
    Therefore, by Lemma~\ref{lem:conditional-many}, both~$f(x)\mid f(x')$ and~$f(x)\mid f(x'),f(x'')$ are (univariate) Gaussian.
    Choosing~$z_1=f(x)$ and~$z_2=f(x')$ in the statement of Lemma~\ref{lem:conditional-many} yields~\eqref{eq:conditional-mean-1} and~\eqref{eq:conditional-variance-1}, while choosing~$z_2=(f(x'),f(x''))$ yields
    \begin{align*}
        \E[f(x)\mid f(x'),f(x')]
        &= m(x)+\begin{bmatrix} C(x,x') \\ C(x,x'') \end{bmatrix}^T
            \begin{bmatrix} V(x') & C(x',x'') \\ C(x',x'') & V(x'') \end{bmatrix}^{-1}
            \begin{bmatrix} f(x')-m(x') \\ f(x'')-m(x'') \end{bmatrix} \\
        &= m(x)+w(x,x',x'')^T\begin{bmatrix} f(x')-m(x') \\ f(x'')-m(x'') \end{bmatrix}
    \end{align*}
    and
    \begin{align*}
        \Var(f(x)\mid f(x'),f(x'))
        &= V(x)-\begin{bmatrix} C(x,x') \\ C(x,x'') \end{bmatrix}^T
            \begin{bmatrix} V(x') & C(x',x'') \\ C(x',x'') & V(x'') \end{bmatrix}^{-1}
            \begin{bmatrix} C(x,x') \\ C(x,x'') \end{bmatrix} \\
        &= V(x)-w(x,x',x'')^T\begin{bmatrix} C(x,x') \\ C(x,x'') \end{bmatrix}.\qedhere
    \end{align*}
\end{proof}
\fi

For example, suppose~$x_k\le x\le x_{k+1}$ for some~$k<n$.
Lemma~\ref{lem:conditional-12} characterizes the distributions of~$f(x)\mid f(x_{k+1})$ and~$f(x)\mid f(x_k),f(x_{k+1})$ when the values of~$f(x_k)$ and~$f(x_{k+1})$ are known.
However, variation in the errors~$\epsilon_i$ makes the values of~$f(x_k)$ and~$f(x_{k+1})$ unknown.
So the conditional means of~$f(x)\mid f(x_{k+1})$ and~$f(x)\mid f(x_k),f(x_{k+1})$ given~$\Dcal$ are random.
But their conditional variances given~$\Dcal$ are \emph{not} random; by Lemma~\ref{lem:conditional-12}, these variances depend on only the known values of the variance function~$V:[0,\infty)\to\R$ and covariance function~$C:[0,\infty)\to\R$.

\section{Conditional moments of a Gauss-Markov process}
\label{sec:arbitrary}

Theorem~\ref{thm:moments-arbitrary} refines Lemma~\ref{lem:conditional-many} by imposing properties~(G) and~(M).
Specifically, it characterizes the mean and variance of~$f(x)\mid\Dcal$ when~$\{f(x)\}_{x\ge0}$ is a sample path of an arbitrary Gauss-Markov process.
These moments depend on the location of~$x$ relative to the points~$x_i$ at which~$\Dcal$ contains noisy observations~$y_i$ of~$f(x_i)$.

\begin{theorem}
    \label{thm:moments-arbitrary}
    Let~$\{f(x)\}_{x\ge0}$ be a sample path of a Gauss-Markov process.
    Suppose~$\Dcal=\{(x_i,y_i)\}_{i=1}^n$ is generated by the process~\eqref{eq:dgp}, where~$x_i\ge0$ is non-decreasing in~$i$ and the errors~$\epsilon_i=y_i-f(x_i)$ are jointly Gaussian with~$\E[\epsilon_i\mid x_1,x_2\ldots,x_n]=0$ independently of~$\{f(x)\}_{x\ge0}$.
    Then~$f(x)\mid\Dcal$ is Gaussian for all~$x\ge0$.
    Moreover:
    \begin{enumerate}

        \item[(i)]
        If~$x\le x_1$, then
        \[ \E[f(x)\mid\Dcal]=m(x)+\frac{C(x_1,x)}{V(x_1)}(\E[f(x_1)\mid\Dcal]-m(x_1)) \]
        and
        \[ \Var(f(x)\mid\Dcal)=\left(\frac{C(x_1,x)}{V(x_1)}\right)^2\Var(f(x_1)\mid\Dcal)+\Var\left(f(x)\mid f(x_1)\right). \]

        \item[(ii)]
        If~$x_k\le x\le x_{k+1}$ for some~$k<n$, then
        \[ \E[f(x)\mid\Dcal]=m(x)+w(x,x_k,x_{k+1})^T\begin{bmatrix} \E[f(x_k)\mid\Dcal]-m(x_k) \\ \E[f(x_{k+1})\mid\Dcal]-m(x_{k+1})\end{bmatrix} \]
        and
        \begin{align*}
            \Var(f(x)\mid\Dcal)
            &= w(x,x_k,x_{k+1})^T\Var\left(\begin{bmatrix} f(x_k) \\ f(x_{k+1}) \end{bmatrix}\mid\Dcal\right)w(x,x_k,x_{k+1}) \\
            &\quad +\Var\left(f(x)\mid f(x_k),f(x_{k+1})\right)
        \end{align*}
        with~$w(x,x_k,x_{k+1})$ defined as in Lemma~\ref{lem:conditional-12}.

        \item[(iii)]
        If~$x\ge x_n$, then
        \[ \E[f(x)\mid\Dcal]=m(x)+\frac{C(x_n,x)}{V(x_n)}(\E[f(x_n)\mid\Dcal]-m(x_n)) \]
        and
        \[ \Var(f(x)\mid\Dcal)=\left(\frac{C(x_n,x)}{V(x_n)}\right)^2\Var(f(x_n)\mid\Dcal)+\Var\left(f(x)\mid f(x_n)\right). \]

    \end{enumerate}
\end{theorem}
\ifbodyproofs\begin{proof}[Proof of Theorem~\ref{thm:moments-arbitrary}]
    The~$y_i=f(x_i)+\epsilon_i$ are sums of (jointly) Gaussian random variables, so they are also jointly Gaussian.
    Thus~$(y_1,y_2,\ldots,y_n,f(x))$ is multivariate Gaussian.
    It follows from Lemma~\ref{lem:conditional-many} that~$f(x)\mid\Dcal=f(x)\mid(y_1,y_2,\ldots,y_n)$ is (univariate) Gaussian.

    By property~(M) and the errors' independence of~$\{f(x)\}_{x\ge0}$, we know~$f(x)$ is conditionally independent of~$\Dcal$ given the values~$f(x_i)$ with~$x_i$ closest to~$x$.
    I use this fact to prove cases~(i)--(iii):
    \begin{enumerate}

        \item[(i)]
        Suppose~$x\le x_1$.
        Then~$f(x)$ is conditionally independent of~$\Dcal$ given~$f(x_1)$.
        Thus
        \begin{align*}
            \E[f(x)\mid\Dcal]
            &\overset{\star}{=} \E\left[\E[f(x)\mid\Dcal_n,f(x_1)]\mid\Dcal\right] \\
            &= \E\left[\E[f(x)\mid f(x_1)]\mid\Dcal\right] \\
            &\overset{\star\star}{=} \E\left[m(x)+\frac{C(x_1,x)}{V(x_1)}(f(x_1)-m(x_1))\mid\Dcal\right] \\
            &= m(x)+\frac{C(x_1,x)}{V(x_1)}(\E[f(x_1)\mid\Dcal]-m(x_1)),
        \end{align*}
        where~$\star$ holds by the law of total expectation and~$\star\star$ holds by Lemma~\ref{lem:conditional-12}.
        Likewise
        \begin{align*}
            \Var(f(x)\mid\Dcal)
            &\overset{\star}{=} \Var\left(\E[f(x)\mid\Dcal,f(x_1)]\mid\Dcal\right)+\E\left[\Var(f(x)\mid\Dcal,f(x_1))\mid\Dcal\right] \\
            &= \Var\left(\E[f(x)\mid f(x_1)]\mid\Dcal\right)+\E\left[\Var(f(x)\mid f(x_1))\mid\Dcal\right] \\
            &\overset{\star\star}{=} \Var\left(m(x)+\frac{C(x_1,x)}{V(x_1)}(f(x_1)-m(x_1))\mid\Dcal\right)+\Var(f(x)\mid f(x_1)) \\
            &= \left(\frac{C(x_1,x)}{V(x_1)}\right)^2\Var(f(x_1)\mid\Dcal)+\Var(f(x)\mid f(x_1)),
        \end{align*}
        where~$\star$ holds by the law of total variance and~$\star\star$ holds by Lemma~\ref{lem:conditional-12}.

        \item[(ii)]
        Now suppose~$x_k\le x\le x_{k+1}$ for some~$k<n$.
        Then~$f(x)$ is conditionally independent of~$\Dcal$ given~$f(x_k)$ and~$f(x_{k+1})$.
        Thus
        \begin{align*}
            \E[f(x)\mid\Dcal]
            &\overset{\star}{=} \E\left[\E[f(x)\mid\Dcal_n,f(x_k),f(x_{k+1})]\mid\Dcal\right] \\
            &= \E\left[\E[f(x)\mid f(x_k),f(x_{k+1})]\mid\Dcal\right] \\
            &\overset{\star\star}{=} \E\left[m(x)+w(x,x_k,x_{k+1})^T\begin{bmatrix} f(x_k)-m(x_k) \\ f(x_{k+1})-m(x_{k+1}) \end{bmatrix}\mid\Dcal\right] \\
            &= m(x)+w(x,x_k,x_{k+1})^T\begin{bmatrix} \E[f(x_k)\mid\Dcal]-m(x_k) \\ \E[f(x_{k+1})\mid\Dcal]-m(x_{k+1}) \end{bmatrix},
        \end{align*}
        where~$\star$ holds by the law of total expectation and~$\star\star$ holds by Lemma~\ref{lem:conditional-12}.
        Likewise
        \begin{align*}
            \Var(f(x)\mid\Dcal)
            &\overset{\star}{=} \Var\left(\E\left[f(x)\mid\Dcal,f(x_k),f(x_{k+1})\right]\mid\Dcal\right)+\E\left[\Var\left(f(x)\mid\Dcal,f(x_k),f(x_{k+1})\right)\mid\Dcal\right] \\
            &= \Var\left(\E\left[f(x)\mid f(x_k),f(x_{k+1})\right]\mid\Dcal\right)+\Var\left(f(x)\mid f(x_k),f(x_{k+1})\right) \\
            &\overset{\star\star}{=} \Var\left(m(x)+w(x,x',x'')^T\begin{bmatrix} f(x')-m(x') \\ f(x'')-m(x'') \end{bmatrix}\mid\Dcal\right) \\
            &\quad +\Var\left(f(x)\mid f(x_k),f(x_{k+1})\right) \\
            &= w(x,x_k,x_{k+1})^T\Var\left(\begin{bmatrix} f(x_k) \\ f(x_{k+1}) \end{bmatrix}\mid\Dcal\right)w(x,x_k,x_{k+1}) \\
            &\quad +\Var\left(f(x)\mid f(x_k),f(x_{k+1})\right),
        \end{align*}
        where~$\star$ holds by the law of total variance and~$\star\star$ holds by Lemma~\ref{lem:conditional-12}.

        \item[(iii)]
        Finally, suppose~$x\ge x_n$.
        Then~$f(x)$ is conditionally independent of~$\Dcal$ given~$f(x_n)$.
        The result follows from similar arguments used to prove case~(i).\qedhere

    \end{enumerate}
\end{proof}
\fi

If~$x\le x_1$, then the estimate~$\hat{f}_\Dcal(x)\equiv\E[f(x)\mid\Dcal]$ of~$f(x)$ is linear in the estimate~$\hat{f}_\Dcal(x_1)$ of~$f(x_1)$.
Conversely, if~$x\ge x_n$, then~$\hat{f}_\Dcal(x)$ is linear in the estimate~$\hat{f}_\Dcal(x_n)$ of~$f(x_n)$.
In both cases, we can construct~$\hat{f}_\Dcal(x)$ in two steps:
\begin{enumerate}

    \item
    Use~$\Dcal$ to estimate the ``boundary'' values~$f(x_1)$ and~$f(x_n)$;

    \item
    Extrapolate~$\hat{f}_\Dcal(x)$ from the closest boundary estimate.

\end{enumerate}

For example, suppose~$\Dcal=\{(x_1,y_1)\}$ contains one observation~$y_1=f(x_1)+\epsilon_1$ with error~$\epsilon_1\sim\Ncal(0,\veps)$.
Then
\[ f(x_1)\mid\Dcal_1\sim\Ncal\left(m(x_1)+\frac{V(x_1)}{V(x_1)+\veps}(y_1-m(x_1)),\ \left(\frac{1}{V(x_1)}+\frac{1}{\veps}\right)^{-1}\right) \]
by Lemma~\ref{lem:conditional-many}, and so
\begin{align*}
    \hat{f}_\Dcal(x)
    &= m(x)+\frac{C(x_1,x)}{V(x_1)}\left(\E[f(x_1)\mid\Dcal]-m(x_1)\right) \\
    &= m(x)+\frac{C(x_1,x)}{V(x_1)+\veps}(y_1-m(x_1))
\end{align*}
is linear in the deviation~$(y_1-m(x_1))$ of~$y_1$ from its mean~$m(x_1)$.
Intuitively, this deviation provides information about the difference~$(f(x_1)-m(x_1))$, which the factor~$C(x_1,x)/(V(x)+\veps)$ translates into information about the difference~$(f(x)-m(x))$.
This factor is larger when the covariance~$C(x_1,x)$ of~$y_1$ and~$f(x)$ is larger, and when the variance~$V(x)+\veps$ of~$y_1$ is smaller.

If~$n>1$ and~$x_1\le x\le x_n$, then we can construct~$\hat{f}_\Dcal(x)$ in three steps:
\begin{enumerate}

    \item
    Find an index~$k<n$ for which~$x_k\le x\le x_{k+1}$;

    \item
    Use~$\Dcal$ to estimate the values of~$f(x_k)$ and~$f(x_{k+1})$;

    \item
    Interpolate~$\hat{f}_\Dcal(x)$ from the estimates~$\hat{f}_\Dcal(x_k)$ and~$\hat{f}_\Dcal(x_{k+1})$.

\end{enumerate}
The interpolant~$\hat{f}_\Dcal(x)$ is a weighted sum of the mean~$m(x)$, the deviation~$(\hat{f}_\Dcal(x_k)-m(x_k))$, and the deviation~$(\hat{f}_\Dcal(x_{k+1})-m(x_{k+1}))$.
The weights on the two deviations depend on the (co)variances of~$f(x)$, $f(x_k)$, and~$f(x_{k+1})$, as well as the (co)variances of~$\epsilon_k$ and~$\epsilon_{k+1}$.

Theorem~\ref{thm:moments-arbitrary} expresses the moments of~$f(x)\mid\Dcal$ in terms of the moments of the~$f(x_i)\mid\Dcal$.
The latter moments can be computed analytically if~$n$ is small (e.g., as in the case with~$n=1$ above), or numerically if~$n$ is large.
The expressions in Theorem~\ref{thm:moments-arbitrary} reveal how changing the moments of the~$f(x_i)\mid\Dcal$ changes the moments of~$f(x)\mid\Dcal$.

\section{Conditional moments of a Brownian motion}
\label{sec:brownian}

I now consider a specific Gauss-Markov process: a Brownian motion with known drift~$\mu\in\R$ and scale~$\sigma\ge0$, and unknown initial value~$f(0)\sim\Ncal(\mu_0,\vz)$.%
\footnote{
So the sample path~$\{f(x)\}_{x\ge0}$ solves the stochastic differential equation
\[ \der f(x)=\mu\,\der x+\sigma\,\der W(x), \]
where~$\{W(x)\}_{x\ge0}$ is an unknown sample path of a (standard) Wiener process.
This process has initial value~$W(0)=0$ and iid Gaussian increments~$\der W(x)\equiv W(x+\der x)-W(x)\sim\Ncal(0,\der x)$.
}
This process has mean and covariance functions defined by
\[ m(x)=\mu_0+\mu x \]
and
\[ C(x,x')=\vz+\vm\min\{x,x'\} \]
for all~$x,x'\ge0$.%
\footnote{
Thus~$V(x)=\vz+\vm x$ for all~$x\ge0$.
}
Substituting these expressions into Theorem~\ref{thm:moments-arbitrary} yields the following result.

\begin{corollary}
    \label{crly:moments-brownian}
    Let~$\{f(x)\}_{x\ge0}$ be a sample path of a Brownian motion with drift~$\mu\in\R$, scale~$\sigma\ge0$, and initial value~$f(0)\sim\Ncal(\mu_0,\vz)$.
    Define~$\Dcal=\{(x_i,y_i)\}_{i=1}^n$ as in Theorem~\ref{thm:moments-arbitrary}.
    Then~$f(x)\mid\Dcal$ is Gaussian for all~$x\ge0$.
    Moreover:
    \begin{enumerate}

        \item[(i)]
        If~$x\le x_1$, then
        \[ \E[f(x)\mid\Dcal]=\mu_0+\mu x+\frac{\vz+\vm x}{\vz+\vm x_1}(\E[f(x_1)\mid\Dcal]-\mu_0-\mu x_1) \]
        and
        \[ \Var(f(x)\mid\Dcal)=\left(\frac{\vz+\vm x}{\vz+\vm x_1}\right)^2\Var(f(x_1)\mid\Dcal)+\left(\frac{\vz+\vm x}{\vz+\vm x_1}\right)\vm(x_1-x). \]

        \item[(ii)]
        If~$x_k\le x\le x_{k+1}$ for some~$k<n$, then
        \[ \E[f(x)\mid\Dcal]=\frac{x_{k+1}-x}{x_{k+1}-x_k}\E[f(x_k)\mid\Dcal]+\frac{x-x_k}{x_{k+1}-x_k}\E[f(x_{k+1})\mid\Dcal] \]
        and
        \begin{align*}
            \Var(f(x)\mid\Dcal)
            &= \frac{1}{(x_{k+1}-x_k)^2}
                \begin{bmatrix} x_{k+1}-x \\ x-x_k \end{bmatrix}^T
                \Var\left(\begin{bmatrix} f(x_k) \\ f(x_{k+1}) \end{bmatrix}\mid\Dcal\right)
                \begin{bmatrix} x_{k+1}-x \\ x-x_k \end{bmatrix} \\
            &\quad +\frac{\vm(x_{k+1}-x)(x-x_k)}{x_{k+1}-x_k}.
        \end{align*}

        \item[(iii)]
        If~$x\ge x_n$, then
        \[ \E[f(x)\mid\Dcal]=\E[f(x_n)\mid\Dcal]+\mu(x-x_n) \]
        and
        \[ \Var(f(x)\mid\Dcal)=\Var(f(x_n)\mid\Dcal)+\vm(x-x_n). \]

    \end{enumerate}
\end{corollary}
\ifbodyproofs\begin{proof}[Proof of Corollary~\ref{crly:moments-brownian}]
    The Brownian motion is a Gauss-Markov process, so~$f(x)\mid\Dcal$ is Gaussian by Theorem~\ref{thm:moments-arbitrary}.
    We can also use that theorem to prove cases~(i)--(iii).

    Consider cases~(i) and~(iii).
    By Lemma~\ref{lem:conditional-12} and the definitions of~$m$ and~$C$, we have
    %
    \begin{align}
        \Var(f(x)\mid f(x'))
        &= \vz+\vm x-\frac{(\vz+\vm\min\{x,x'\})^2}{\vz+\vm x'} \notag \\
        &= \frac{\vz\vm\left(x+x'-2\min\{x,x'\}\right)+\sm^4\left(xx'-\min\{x,x'\}^2\right)}{\vz+\vm x'} \notag \\
        &= \frac{\vz\vm\abs{x-x'}+\sm^4\min\{x,x'\}\abs{x-x'}}{\vz+\vm x'} \notag \\
        &= \left(\frac{\vz+\vm\min\{x,x'\}}{\vz+\vm x'}\right)\vm\abs{x-x'} \notag \\
        &= \vm\abs{x-x'}\begin{cases}
            \frac{\vz+\vm x}{\vz+\vm x'} & \text{if}\ x\le x' \\
            1 & \text{otherwise}
        \end{cases} \label{eq:conditional-variance-1-brownian}
    \end{align}
    for all~$x'\ge0$.
    Cases~(i) and~(iii) follow from substituting~\eqref{eq:conditional-variance-1-brownian}, along the definitions of~$m$ and~$C$, into the statement of Theorem~\ref{thm:moments-arbitrary}.

    Now consider case~(ii).
    Define~$w(x,x_k,x_{k+1})$ as in Lemma~\ref{lem:conditional-12}.
    If~$x_k\le x\le x_{k+1}$, then
    \begin{align*}
        w(x,x_k,x_{k+1})
        &\equiv \frac{1}{V(x_k)V(x_{k+1})-C(x_k,x_{k+1})^2}
            \begin{bmatrix} C(x,x_k)V(x_{k+1})-C(x,x_{k+1})C(x_k,x_{k+1}) \\ C(x,x_{k+1})V(x_k)-C(x,x_k)C(x_k,x_{k+1}) \end{bmatrix} \\
        &= \frac{1}{V(x_k)V(x_{k+1})-V(x_k)^2}
            \begin{bmatrix} V(x_k)V(x_{k+1})-V(x)V(x_k) \\ V(x)V(x_k)-V(x_k)V(x_k) \end{bmatrix} \\
        &= \frac{1}{x_{k+1}-x_k}\begin{bmatrix} x_{k+1}-x \\ x-x_k \end{bmatrix}
    \end{align*}
    because~$C(x',x'')=V(\min\{x',x''\})$ for all~$x',x''\ge0$ by the definitions of~$C$ and~$V$.
    So
    \begin{align*}
        &w^T\begin{bmatrix} \E[f(x_k)\mid\Dcal]-m(x_k) \\ \E[f(x_{k+1})\mid\Dcal]-m(x_{k+1})\end{bmatrix} \\
        &= \frac{1}{x_{k+1}-x_k}
            \begin{bmatrix} x_{k+1}-x \\ x-x_k \end{bmatrix}^T
            \begin{bmatrix} \E[f(x_k)\mid\Dcal] \\ \E[f(x_{k+1})\mid\Dcal]\end{bmatrix}
        -\frac{1}{x_{k+1}-x_k}
            \begin{bmatrix} x_{k+1}-x \\ x-x_k \end{bmatrix}^T
            \begin{bmatrix} \mu_0+\mu x_k \\ \mu_0+\mu x_{k+1} \end{bmatrix} \\
        &= \frac{(x_{k+1}-x)\E[f(x_k)\mid\Dcal]+(x-x_k)\E[f(x_{k+1})\mid\Dcal]}{x_{k+1}-x_k}
            -\mu_0-\mu x
    \end{align*}
    by the definition of~$m$, and
    \begin{align*}
        \Var(f(x)\mid f(x_k),f(x_{k+1}))
        &= V(x)-w(x,x_k,x_{k+1})^T\begin{bmatrix} C(x,x_k) \\ C(x,x_{k+1}) \end{bmatrix} \\
        &= \vz+\vm x-\frac{1}{x_{k+1}-x_k}\begin{bmatrix} x_{k+1}-x \\ x-x_k \end{bmatrix}^T\begin{bmatrix} \vz+\vm x_k \\ \vz+\vm x \end{bmatrix} \\
        &= \frac{\vm(x_{k+1}-x)(x-x_k)}{x_{k+1}-x_k}
    \end{align*}
    by Lemma~\ref{lem:conditional-12}.
    The result follows from Theorem~\ref{thm:moments-arbitrary}(ii).
\end{proof}
\fi

If~$\{f(x)\}_{x\ge0}$ is a sample path of a Brownian motion, then the estimate~$\hat{f}_\Dcal(x)\equiv\E[f(x)\mid\Dcal]$ of~$f(x)$ is piecewise linear in~$x\ge0$.
For example, consider the canonical case in which the initial value~$f(0)=\mu_0$ is known.
As~$x$ increases from zero, the estimate
\[ \hat{f}_\Dcal(x)=\begin{cases}
    \mu_0+\frac{x}{x_1}\left(\hat{f}_\Dcal(x_1)-\mu_0\right) & \text{if}\ x<x_1 \\
    \hat{f}_\Dcal(x_k)+\frac{x-x_k}{x_{k+1}-x_k}\left(\hat{f}_\Dcal(x_{k+1})-\hat{f}_\Dcal(x_k)\right) & \text{if}\ x_k\le x<x_{k+1}\ \text{for some}\ k<n \\
    \hat{f}_\Dcal(x_n)+\mu(x-x_n) & \text{otherwise}
    \end{cases} \]
interpolates linearly between~$(0,\mu_0)$ and~$(x_1,\hat{f}_\Dcal(x_1))$, then between~$(x_1,\hat{f}_\Dcal(x_1))$ and~$(x_2,\hat{f}_\Dcal(x_2))$, then between~$(x_2,\hat{f}_\Dcal(x_2))$ and~$(x_3,\hat{f}_\Dcal(x_3))$, and so on until~$(x_n,\hat{f}_\Dcal(x_n))$.
Beyond this point, the data~$\Dcal$ provide no information about~$f(x)$.
Consequently, the estimate~$\hat{f}_\Dcal(x)$ treats~$\{f(x)\}_{x\ge x_n}$ as a Brownian motion with drift~$\mu$, scale~$\sm$, and (possibly random) initial value~$\hat{f}_\Dcal(x_n)$.

\begin{figure}
    \centering
    \includegraphics[width=0.75\linewidth]{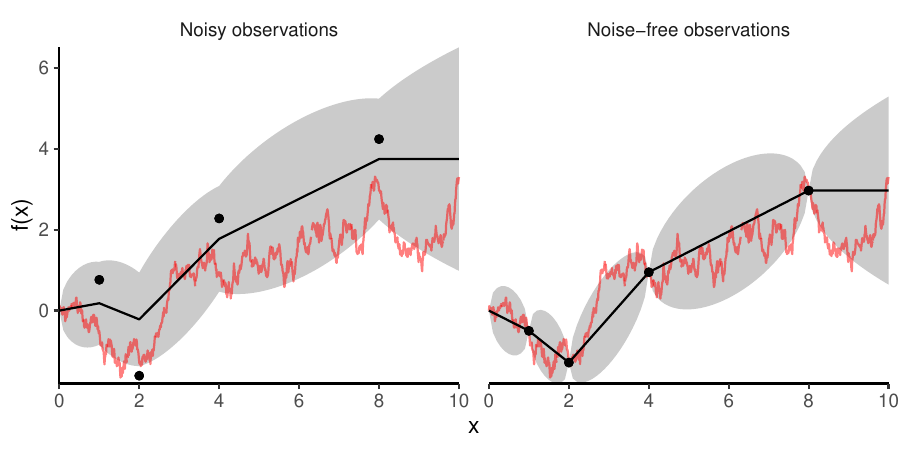}
    \caption{Learning about a Brownian motion from noisy and noise-free observations}
    \label{fig:brownian}
    \caption*{\footnotesize{\itshape Notes}:
    Red lines show sample path~$\{f(x)\}_{x\ge0}$ of Brownian motion with drift~$\mu=0$, scale~$\sm=1$, and initial value~$f(0)=0$.
    Black dots represent observations~$(x_i,y_i)$ with~$x_i=2^{i-1}$ and $y_i=f(x_i)+\epsilon_i$.
    Black lines represent estimates~$\hat{f}_\Dcal(x)$ of~$f(x)$ when~$\Dcal=\{(x_i,y_i)\}_{i=1}^4$.
    Gray regions represent 90\% confidence intervals, constructed analytically from the conditional variances defined in Corollary~\ref{crly:moments-brownian}.
    Left panel has~$\epsilon_i\sim\Ncal(0,1)$ iid; right panel has~$\epsilon_i=0$ for each~$i$.
    }
\end{figure}

I illustrate this behavior in Figure~\ref{fig:brownian}.
It shows how~$\hat{f}_\Dcal(x)$ varies with~$x\ge0$ when~$(\mu_0,\mu,\sz,\sm)=(0,0,0,1)$ and the data contain~$n=4$ observations with iid errors.
Figure~\ref{fig:brownian} also shows the 90\% confidence interval around~$\hat{f}_\Dcal(x)$, constructed analytically from the conditional variances defined in Corollary~\ref{crly:moments-brownian}.
This interval expands as~$x$ moves away from the sampled points~$x_i$.%
\footnote{
If~$f(0)$ is known (i.e., $\vz=0$) and the observations~$y_i$ have no noise (i.e., $\epsilon_i=0$ for each~$i$), then the MSE
\[ \E\left[\left(f(x)-\hat{f}_\Dcal(x)\right)^2\mid\Dcal\right]=\vm\begin{cases}
    \frac{x(x_1-x)}{x_1} & \text{if}\ x<x_1 \\
    \frac{(x_{k+1}-x)(x-x_k)}{x_{k+1}-x_k} & \text{if}\ x_k\le x<x_{k+1}\ \text{for some}\ k<n \\
    x-x_n & \text{otherwise}
    \end{cases} \]
attains its piecewise maxima at the midpoint of each piece.
}

Suppose the data are \emph{not} noisy (i.e., $\epsilon_i=0$ for each~$i$) and~$x_k\le x\le x_{k+1}$ for some~$k<n$.
Then the distribution of~$f(x)\mid\Dcal$ coincides with the distribution obtained by assuming~$\{f(x)\}_{x_k\le x\le x_{k+1}}$ is a sample path of a Brownian bridge with scale~$\sm$, known initial value~$f(x_k)$, and known terminal value~$f(x_{k+1})$.%
\footnote{
See, e.g., Section~5.6.B of \cite{Karatzas-Shreve-1988-} for more information about Brownian bridges.
}
Corollary~\ref{crly:moments-brownian}(ii) generalizes to Brownian bridges with \emph{unknown} initial and terminal values.
For example, suppose~$x_1\le x\le x_2$ and that
\[ \begin{bmatrix} f(x_1) \\ f(x_2) \end{bmatrix}\mid\Dcal\sim\Ncal\left(
    \begin{bmatrix} \mu_1 \\ \mu_2 \end{bmatrix},\, 
    \begin{bmatrix}
        \sigma_1^2 & \rho\sigma_1\sigma_2 \\
        \rho\sigma_1\sigma_2 & \sigma_2^2
    \end{bmatrix}
    \right) \]
for some means~$\mu_1,\mu_2\in\R$, variances~$\sigma_1^2,\sigma_2^2\ge0$, and correlation~$\rho\in[-1,1]$.
Then~$f(x)\mid\Dcal$ has mean
\[ \E[f(x)\mid\Dcal]=\mu_1+\frac{x-x_1}{x_2-x_1}\left(\mu_2-\mu_1\right) \]
and variance
\begin{align*}
    \Var(f(x)\mid\Dcal)
    &= \frac{(1-x)^2\sigma_1^2+2x(1-x)\rho\sigma_1\sigma_2+x^2\sigma_2^2}{(x_2-x_1)^2}+\frac{\vm(x_2-x)(x-x_1)}{x_2-x_1}.
\end{align*}
Choosing~$\sigma_1=0$ and~$\sigma_2=0$ yields the mean and variance for the Brownian bridge on~$[x_1,x_2]$ with~$f(x_1)=\mu_1$ and~$f(x_2)=\mu_2$.

{%
\small
\raggedright
\bibliographystyle{apalike}
\bibliography{references}
}

\appendix

\clearpage
\section{Proofs}
\label{sec:proofs}

\linespread{1}

\ifbodyproofs\else\begin{proof}[Proof of Lemma~\ref{lem:conditional-12}]
    Now~$f(x)$, $f(x')$, and~$f(x'')$ are jointly Gaussian by property~(G).
    Therefore, by Lemma~\ref{lem:conditional-many}, both~$f(x)\mid f(x')$ and~$f(x)\mid f(x'),f(x'')$ are (univariate) Gaussian.
    Choosing~$z_1=f(x)$ and~$z_2=f(x')$ in the statement of Lemma~\ref{lem:conditional-many} yields~\eqref{eq:conditional-mean-1} and~\eqref{eq:conditional-variance-1}, while choosing~$z_2=(f(x'),f(x''))$ yields
    \begin{align*}
        \E[f(x)\mid f(x'),f(x')]
        &= m(x)+\begin{bmatrix} C(x,x') \\ C(x,x'') \end{bmatrix}^T
            \begin{bmatrix} V(x') & C(x',x'') \\ C(x',x'') & V(x'') \end{bmatrix}^{-1}
            \begin{bmatrix} f(x')-m(x') \\ f(x'')-m(x'') \end{bmatrix} \\
        &= m(x)+w(x,x',x'')^T\begin{bmatrix} f(x')-m(x') \\ f(x'')-m(x'') \end{bmatrix}
    \end{align*}
    and
    \begin{align*}
        \Var(f(x)\mid f(x'),f(x'))
        &= V(x)-\begin{bmatrix} C(x,x') \\ C(x,x'') \end{bmatrix}^T
            \begin{bmatrix} V(x') & C(x',x'') \\ C(x',x'') & V(x'') \end{bmatrix}^{-1}
            \begin{bmatrix} C(x,x') \\ C(x,x'') \end{bmatrix} \\
        &= V(x)-w(x,x',x'')^T\begin{bmatrix} C(x,x') \\ C(x,x'') \end{bmatrix}.\qedhere
    \end{align*}
\end{proof}
\fi

\ifbodyproofs\else\fi

\ifbodyproofs\else\fi

\end{document}